\documentclass{birkjour_t2}
\usepackage{amssymb}
\usepackage{amsthm}
\usepackage{bm}
\usepackage{amsmath}
\usepackage{enumerate}
\usepackage{enumitem}
 \usepackage{lineno}

 \theoremstyle{definition}
 
 \theoremstyle{remark}

 \numberwithin{equation}{section}

\begin{document}

\title[Stress-rate type model of strain-limiting viscoelasticity]{A thermodynamically consistent stress-rate type model of one-dimensional strain-limiting viscoelasticity}
\author[H. A. Erbay]{H. A. Erbay}
    \address{Department of Natural and Mathematical Sciences \\ Faculty of Engineering \\ Ozyegin University\\ Cekmekoy 34794 Istanbul\\ Turkey \\
    {ORCID ID: 0000-0002-5167-609X}}
    \email{husnuata.erbay@ozyegin.edu.tr}
    
\author[Y. \c Seng\"{u}l]{Y. \c Seng\"{u}l${^{*}}$}
    \address{Faculty of Engineering and Natural Sciences\\ Sabanci University\\ Tuzla 34956 Istanbul \\  Turkey\\ 
    {ORCID ID: 0000-0001-5923-3173}}
     \email{yaseminsengul@sabanciuniv.edu}
     \thanks{${^{*}}$ corresponding author}


\begin{abstract}
 We  introduce a one-dimensional stress-rate type nonlinear viscoelastic model for solids that obey the assumptions of the strain-limiting theory. Unlike the classical   viscoelasticity  theory, the critical hypothesis in the present strain-limiting theory is that the  linearized strain  depends nonlinearly on the stress and the stress rate.  We  show the thermodynamic consistency of the model using the complementary free energy, or equivalently, the Gibbs free energy. This allows us to take the stress and the stress rate as primitive variables instead of kinematical quantities such as deformation or strain. We also show that the non-dissipative part of the materials in consideration have a stored energy.  We  compare the new stress-rate type model with the strain-rate type viscoelastic model due to Rajagopal from the points of view of energy decay, the nonlinear differential equations of motion and Fourier analysis of the corresponding linear models.
\end{abstract}

\maketitle


\section{Introduction}\label{sec:sec1}

In the present work, our main aim is to introduce a thermodynamically consistent stress-rate model for the one-dimensional nonlinear response of strain-limiting viscoelastic solids, and compare it from different perspectives with the strain-rate model of strain-limiting viscoelastic solids.

In the classical (Cauchy) theory of elasticity, it is postulated that the stress depends (possibly nonlinearly) on the deformation gradient or on  an appropriate strain measure. In a series of articles \cite{Rajagopal2003,Rajagopal2007a,Rajagopal2010a,Rajagopal2011a,Rajagopal2011b,Rajagopal2014a,Rajagopal2018}, a new model of elastic materials (hereafter called  the strain-limiting model of elasticity) was developed by Rajagopal, in which that the linearized strain is a nonlinear function of the stress. That is, in this model, the displacement gradients are assumed to be small, while placing no restriction on the stress. The motivation for introducing the strain-limiting model of elasticity was two-fold. Firstly, there are experimental observations reported for elastic materials in which a nonlinear relationship between the linearized strain and the stress holds. As a typical example, one can consider the cases at  crack tips in brittle elastic bodies, where the stresses can be large but the strains remain small. Secondly, it is more appropriate to use the stress as a primitive variable instead of a kinematical quantity (deformation). The reason that deformation was not chosen as a primitive variable is that the stress, which is a consequence of  the contact force, causes deformation, that is, the stress is the cause of deformation but not the effect. For a more detailed overview of the strain-limiting model of elasticity  we refer the reader to   \cite{Rajagopal2003,Rajagopal2007a,Rajagopal2010a,Rajagopal2011a,Rajagopal2011b,Rajagopal2014a,Rajagopal2018}. The literature contains many studies that have considered the strain-limiting model in analysing various issues,   spanning from metallic alloys to biological fibers. We refer the reader to the recent articles \cite{Huang2017,Freed2016,Fu2019} and  references therein for various application areas of strain-limiting theories, such as fracture mechanics and biomechanics.

Even though analysis from many different points of view has been done for strain-limiting models of elastic solids, there are only a few attempts to model viscoelastic solid materials within the context of strain-limiting theory. In \cite{Muliana2013} an integral-type  strain-limiting model, where the linearized strain is related to some integral of the stress and the stress rate, was proposed. In \cite{Rajagopal2014b} a three-dimensional strain-rate type (Kelvin-Voigt type) strain-limiting model where a linear combination of the linearized strain and the linearized strain-rate is a function of the stress was introduced. Later,  for a particular form of the constitutive relation suggested in \cite{Rajagopal2014b}, traveling waves in the one-dimensional model were analyzed  in   \cite{Erbay2015}. Mathematical issues concerning the strain-limiting viscoelastic models were studied in \cite{Bulicek2012,Itou2019}. For a self-contained overview of strain-limiting viscoelasticity, we refer the reader to \cite{Sengul}.

The question that naturally arises is whether it is possible to develop a thermodynamically consistent stress-rate type model of strain-limiting materials. The aim of this work is  to build a thermodynamically consistent model of stress-rate type strain-limiting viscoelastic materials. Here, the term``thermodynamically consistent" means a model with the constitutive relation satisfying the second law of thermodynamics in the form of the Clausius-Duhem inequality. To avoid the introduction of a strain-like quantity as a primitive variable we use a formulation based on the complementary free energy (or equivalently, the Gibbs free energy) instead of the classical formulations based on the Helmholtz energy function.   The key point in the present analysis is the assumption that the thermodynamic potentials depend on the stress and the stress rate.

The structure of  the paper is as follows.  In Section \ref{sec:sec2}, we briefly review the strain-rate type model of  three dimensional strain-limiting viscoelasticity. In Section \ref{sec:sec3}, we introduce a one-dimensional stress-rate type model of strain-limiting viscoelasticity. In Section \ref{sec:sec4}, starting from the Clausius-Duhem inequality we derive the thermodynamically consistent constitutive equation for the stress-rate type viscoelastic model. In Section \ref{sec:sec5}, we compare the strain-rate type and stress-rate type models of strain-limiting viscoelasticity  in terms of energy decay, the nonlinear differential equations of motion and Fourier analysis of the corresponding linear models.

\section{Preliminaries: Strain-rate type model of strain-limiting viscoelasticity}\label{sec:sec2}

The goal of  this section is to provide a brief introduction of both some basic terminology  and the constitutive relation for the strain-rate type strain-limiting viscoelastic materials discussed in \cite{Rajagopal2014b,Erbay2015}.

Consider a continuum whose motion is expressed by the mapping function $\boldsymbol{\chi}$ defined by $\mathbf{x}=\boldsymbol{\chi}(\mathbf{X},t)$ where $\mathbf{X}$ is the initial (undeformed)  position vector of a particle X, $\mathbf{x}$ is the  current (deformed) position vector of the same particle and $t$ denotes the time. The displacement vector $\mathbf{u}$, the deformation gradient tensor $\mathbf{F}$ and the velocity vector $\mathbf{v}$ of the particle X are defined as $\mathbf{u}=\mathbf{x}-\mathbf{X}$, $\mathbf{F}=\partial \mathbf{x}/\partial \mathbf{X}$ and $\mathbf{v}=\dot{\mathbf{x}}$, respectively, where the superposed dot denotes the material time derivative. The  left Cauchy-Green strain tensor $\mathbf{B}$ (also known as the Finger deformation tensor)  is given by $\mathbf{B}=\mathbf{F}\mathbf{F}^{T}$ (``T" stands for the transpose). We conclude this paragraph with several more definitions  that will be used below. The linearized strain tensor $\boldsymbol{\epsilon}$ and the stretching tensor $\mathbf{D}$ (i.e., the symmetric part of the velocity gradient or the rate of deformation tensor)    are defined through
\begin{equation}\label{linearstrain}
    \boldsymbol{\epsilon} = \frac{1}{2} \Big(\frac{\partial \mathbf{u}}{\partial \mathbf{x}}
        + \big(\frac{\partial \mathbf{u}}{\partial \mathbf{x}}\big)^{T}\Big), ~~~~
    \mathbf{D}= \frac{1}{2} \Big(\frac{\partial \mathbf{v}}{\partial \mathbf{x}}
        + \big(\frac{\partial \mathbf{v}}{\partial \mathbf{x}}\big)^{T}\Big),
\end{equation}
respectively. The Almansi-Hamel strain tensor   $\mathbf{A}$ is related to the linearized strain tensor $\boldsymbol{\epsilon}$ and to the  left Cauchy-Green strain tensor $\mathbf{B}$ as follows:
\begin{equation}\label{almansistrain}
    \mathbf{A}=\boldsymbol{\epsilon}-\frac{1}{2}\big(\frac{\partial \mathbf{u}}{\partial \mathbf{x}}\big)\big(\frac{\partial \mathbf{u}}{\partial \mathbf{x}}\big)^{T}, ~~~~
    \mathbf{A}=\frac{1}{2}\big(\mathbf{I}-\mathbf{B}^{-1}\big),
\end{equation}
where $\mathbf{I}$ is the identity tensor and the symbol $(.)^{-1}$ stands for the inverse. The above overview of kinematical definitions will suffice for our purposes.

In the absence of body forces, the balance of linear momentum and the balance of angular momentum are given in terms of the Cauchy stress tensor $\mathbf{T}(\mathbf{x},t)$ by
\begin{equation}
    \textrm{div}~\mathbf{T}=\rho \ddot{\mathbf{u}} ~~\mbox{and}~~\mathbf{T}^{T}=\mathbf{T},
\end{equation}
respectively, where $\rho $ is the density in the current configuration and the symbol ``$\textrm{div}$" denotes the divergence operator  with respect to $\mathbf{x}$.  The balance of mass is given by $\rho~ \mathrm{det~}\mathbf{F}=\rho_{0}$ where $\rho_{0}$ is the reference density. The above equations  must be augmented with a constitutive relation. The strain-limiting theory of elastic solids is  based on the smallness assumption on  the displacement gradient $\partial \mathbf{u}/\partial \mathbf{x}$ for all points in the body throughout the time interval considered. We require that each component of $\partial \mathbf{u}/\partial \mathbf{x}$ is small and in some appropriate norm we write:  $ \Vert \partial \mathbf{u}/\partial \mathbf{x}  \Vert = {\mathcal O}(\delta)$  with $\delta << 1$. By neglecting the  nonlinear terms with ${\mathcal O}(\delta ^{2})$ in $\mathbf{A}$ we may approximate $\mathbf{A}$ by  $\boldsymbol{\epsilon}$. So the smallness assumption  allows us to replace the Cauchy-Green tensor $\mathbf{B}$ by $\mathbf{I} + 2 \boldsymbol{\epsilon}$ in the constitutive relation of the solid material. Starting from a subclass of the general implicit constitutive relation $\mathcal{F}(\mathbf{T}, \mathbf{B}) = \mathbf{0}$ in terms of the conjugate pair ($\mathbf{T},\mathbf{B}$) and using the above small strain assumption,  Rajagopal \cite{Rajagopal2003,Rajagopal2007a} derived the following  constitutive relation for the strain-limiting model of homogeneous isotropic elastic materials
\begin{displaymath}
    \boldsymbol{\epsilon} = \beta_{0} \mathbf{I} + \beta_{1} \mathbf{T} + \beta_{2} \mathbf{T}^{2},
\end{displaymath}
where $\beta_{i} = \beta_{i}(\rho,I_{1}, I_{2}, I_{3})$, $(i = 0, 1, 2)$ are scalar functions of the density and the three invariants of the stress tensor,  $I_{1} = \mathrm{tr~} \mathbf{T}$, $I_{2} = \frac{1}{2} \mathrm{tr~} \mathbf{T}^{2}$, $I_{3} = \frac{1}{3} \mathrm{tr~} \mathbf{T}^{3}$ (here the symbol $\mathrm{tr}$ denotes the trace of a matrix). We note here that when the stress is zero, deformation does not take place and the definitions of $\mathbf{B}$ and $\boldsymbol{\epsilon}$ imply that $\mathbf{B}=\mathbf{I}$ and $\boldsymbol{\epsilon}=\mathbf{0}$, respectively. This imposes definite restrictions on the forms of the coefficient functions $\beta_{i}$.

Within the context of the strain-limiting approach, however, there are not so many studies about constitutive modeling of viscoelastic materials. Muliana et al. \cite{Muliana2013}  developed a quasi-linear viscoelastic model where the linearized strain is expressed as an integral of a nonlinear measure of the stress. As a subclass of the general implicit constitutive relations of the form $\mathcal{F}(\mathbf{T}, \mathbf{B}, \mathbf{D}) = \mathbf{0}$,  Rajagopal and Saccomandi \cite{Rajagopal2014b}  studied the following strain-rate type model of homogeneous isotropic viscoelastic materials
\begin{displaymath}
    \mathbf{B} + \nu \mathbf{D} = \bar{\beta}_{0} \mathbf{I} + \bar{\beta}_{1} \mathbf{T} + \bar{\beta}_{2} \mathbf{T}^{2},
\end{displaymath}
where $\nu$ is a positive constant and $\bar{\beta}_{i} = \bar{\beta}_{i}(I_{1}, I_{2}, I_{3})$ $(i = 0, 1, 2)$. After making the assumption that $ \Vert \partial \mathbf{u}/\partial \mathbf{x}  \Vert = {\mathcal O}(\delta)$ with $\delta << 1$ in some appropriate norm and approximating $\mathbf{D}$ by $\boldsymbol{\epsilon}_{t}=\partial \boldsymbol{\epsilon}/\partial t$, they proposed the strain-limiting viscoelasticity model  given by
\begin{equation}\label{Raj-Sac-model}
    \boldsymbol{\epsilon} + \nu \boldsymbol{\epsilon}_{t} = \beta_{0} \mathbf{I} + \beta_{1} \mathbf{T} + \beta_{2} \mathbf{T}^{2},
\end{equation}
where $\beta_{i} = \beta_{i}(I_{1}, I_{2}, I_{3})$ $(i = 0, 1, 2)$. Erbay and \c{S}eng\"ul \cite{Erbay2015} considered \eqref{Raj-Sac-model} in one space dimension with a general nonlinear right-hand side and suggested
\begin{equation}\label{Raj-model}
    \epsilon + \nu \epsilon_{t} = g(T),
\end{equation}
where $\epsilon(x,t)$ is the linearized strain (connected to the displacement function $u(x,t)$ by $\epsilon=u_{x}$) and $g(\cdot)$ is a nonlinear function of the Cauchy stress $T(x, t)$ with $g(0)=0$. Using this constitutive relation and the equation of motion in dimensionless variables, they studied  traveling wave solutions   for different forms of $g$ and under the assumption of two constant equilibrium states at infinity.

\section{A stress-rate type strain-limiting model of viscoelastic materials}\label{sec:sec3}

To model  the stress-rate type viscoelastic fluids within the context of implicit constitutive theories,  the relation $\mathcal{F}(\mathbf{T}, \dot{\mathbf{T}}, \partial \mathbf{v}/\partial \mathbf{x}) = \mathbf{0}$ has been considered in \cite{Rajagopal2011c}.  In a similar way, in order to model  the stress-rate type viscoelastic solids in the context of implicit constitutive theories, we consider the  relation $\mathcal{F}(\mathbf{T}, \dot{\mathbf{T}}, \mathbf{B}) = \mathbf{0}$. Furthermore, we will restrict our attention to the case   where the strain is given as a nonlinear function of the stress and its time derivative, namely, $\mathbf{B} = \mathcal{H}(\mathbf{T}, \dot{\mathbf{T}})$. Under the assumption of isotropic materials we have (see Spencer \cite{Spencer1971})
\begin{multline} \label{new-imp-cons}
    \mathbf{B} = \alpha_{0} \mathbf{I} + \alpha_{1} \mathbf{T} + \alpha_{2} \dot{\mathbf{T}} + \alpha_{3} \mathbf{T}^{2}
                + \alpha_{4} \dot{\mathbf{T}}^{2} + \alpha_{5} (\mathbf{T}\dot{\mathbf{T}} + \dot{\mathbf{T}} \mathbf{T})
                + \alpha_{6} (\mathbf{T}^{2} \dot{\mathbf{T}} + \dot{\mathbf{T}} \mathbf{T}^{2}) \\
                + \alpha_{7} (\dot{\mathbf{T}}^{2} \mathbf{T} + \mathbf{T} \dot{\mathbf{T}}^{2}) + \alpha_{8} (\mathbf{T}^{2} \dot{\mathbf{T}}^{2} + \dot{\mathbf{T}}^{2} \mathbf{T}^{2})
\end{multline}
with the scalar functions $\alpha_{i} \,(i = 0, \ldots , 8)$ depending on the invariants
\begin{displaymath}
    \mathrm{tr~} \mathbf{T}, \,\mathrm{tr~} \dot{\mathbf{T}}, \,\mathrm{tr~} \mathbf{T}^{2}, \,\mathrm{tr~} \dot{\mathbf{T}}^{2}, \,\mathrm{tr~} \mathbf{T}^{3}, \,\mathrm{tr~} \dot{\mathbf{T}}^{3}, \,\mathrm{tr~} (\mathbf{T}\dot{\mathbf{T}}), \,\mathrm{tr~} (\mathbf{T}^{2} \dot{\mathbf{T}}), \,\mathrm{tr~} (\dot{\mathbf{T}}^{2} \mathbf{T}), \,\mathrm{tr~}(\mathbf{T}^{2} \dot{\mathbf{T}}^{2}).
\end{displaymath}
Since our aim is to develop a model of stress-rate type strain-limiting viscoelasticity, we are now interested in a linearized version of \eqref{new-imp-cons}. Following Rajagopal and Saccomandi (see (2.3) in  \cite{Rajagopal2014b}) we  assume that
$ \Vert \partial \mathbf{u}/\partial \mathbf{x}  \Vert= {\mathcal O}(\delta)$ and $ \Vert \partial \mathbf{v}/\partial \mathbf{x}  \Vert = {\mathcal O}(\delta)$  with $\delta << 1$. This   allows us to replace $\mathbf{B}$ by $\mathbf{I} + 2 \boldsymbol{\epsilon}$ in the constitutive relation.  Secondly, we assume that the convective terms in the expression of the material time derivative of $\mathbf{T}$ can be neglected, that is, we assume that the products of the gradient of the stress and the velocity are small. As a result, the material time derivative becomes the partial derivative with respect to time. Furthermore, we restrict our attention to  one space dimension. So, at the end, we have a constitutive relation of the form $\epsilon = l(T, T_{t})$, where $\epsilon(x, t)$ is the linearized strain and $T(x, t)$ is the stress.  In order to capture the effect of viscosity while keeping the model as simple as possible we consider the subclass that is linear with respect to $T_{t}$:
\begin{equation}\label{new-model}
    \epsilon = h(T) - \gamma \,T_{t},
\end{equation}
where $h(\cdot)$ is a nonlinear function of the Cauchy stress $T$ with $h(0) = 0$, and $\gamma$ is a constant (see the next section for a discussion on its sign). In the next section,  thermodynamic basis is provided for the above model, that is,  it is shown that, for appropriate choices of the function $h$, (\ref{new-model}) is consistent with the first and second laws of thermodynamics.

\section{Thermodynamic consistency of the stress-rate type model}\label{sec:sec4}

In this section, our aim is to demonstrate that the model \eqref{new-model} is consistent with the first and second laws of thermodynamics. It is important to note that we provide here a thermodynamic framework where the stress rather than a kinematical quantity (such as strain) is a primitive variable. In what follows, we will focus mainly on this aspect of the model. For this reason, and for simplicity of presentation we restrict ourselves to the one-dimensional case.

A thermodynamic basis for the implicit model of elasticity, a subject that is closely related to the discussion in this study,  was provided in \cite{Rajagopal2007b,Rajagopal2009}. In recent studies,   thermodynamically consistent models for the responses of rate-type viscoelastic fluids \cite{Rajagopal2011c} and rate-type thermoviscoelastic solids \cite{Rajagopal2013} have been developed  using a Gibbs-free-energy-based formulation. We now follow a similar approach to show that a constitutive relation  for strain-limiting viscoelastic solids can be derived by stating the Clausis-Duhem inequality in terms of the complementary free energy,  or equivalently, in terms of the Gibbs free energy. Indeed, this is the point of departure from the classical theories where  the Helmholtz free energy  is the most commonly used thermodynamic potential.

\subsection{A thermodynamic framework based on the Almansi-Hamel strain}

In one space dimension $x$, the linearized strain $\epsilon(x,t)$ and the Almansi-Hamel strain $A(x,t)$ are related to the displacement function $u(x,t)$ by $\epsilon=\partial u/\partial x$  and $A=\epsilon -(\partial u/\partial x)^{2}/2$, respectively, due to  \eqref{linearstrain} and  \eqref{almansistrain}. The first law of thermodynamics, i.e. the principle of conservation of energy, takes the form
\begin{equation}\label{cons-energy}
    \rho \dot{\omega} = T \dot{A} - \frac{\partial q}{\partial x},
\end{equation}
where $\rho$ is the density, $T$ is the Cauchy stress, $q$ is the heat flux, $\omega$ is the internal energy per unit mass and all these quantities are functions of $x$ and $t$. In the absence of an entropy source the well-known Clausis-Duhem inequality (see, e.g., \cite{Christensen1971}) states that
\begin{equation}\label{Clausis-Duhem}
    \rho \dot{\eta} + \frac{\partial ~}{\partial x}\left(\frac{q}{\theta}\right) \geq 0,
\end{equation}
where  $\eta$ is the  entropy density and $\theta$ is the absolute temperature. For our purposes, it will be more convenient to express \eqref{Clausis-Duhem} in terms of $\psi$, the Helmhotz free energy  per unit mass,   which is defined by $\psi = \omega - \theta \eta$. Using the time derivative of $\psi = \omega - \theta \eta$   in \eqref{Clausis-Duhem} we get $\rho \dot{\omega} - \rho \dot{\psi} - \rho \eta \dot{\theta} - (q/\theta)  \partial\theta/\partial x + \partial q/\partial x \geq 0$. Furthermore, substituting \eqref{cons-energy} into this inequality gives the Clausis-Duhem inequality in terms of $\psi$ as
\begin{equation}\label{CD-dot}
    -\rho \dot{\psi} + T \dot{A} - \rho \eta \dot{\theta} -  \frac{q}{\theta} \frac{\partial \theta}{\partial x} \geq 0.
\end{equation}
This implies that for constant $A$ and constant $\theta$, the Helmholtz free energy $\psi$ approaches a minimum in equilibrium. If the complementary free energy $\phi_{c}$ (per unit mass)  is defined by
\begin{equation}\label{comp}
    \rho \phi_{c}=TA - \rho \psi,
\end{equation} the Clausis-Duhem inequality \eqref{CD-dot} can be written in terms of $\phi_{c}$ as
\begin{equation}\label{CD-gibbs}
    \rho \dot{\phi_{c}} + TA\frac{\dot{\rho}}{\rho} -\dot{T}A- \rho \eta \dot{\theta}  - \frac{q}{\theta} \frac{\partial \theta}{\partial x} \geq 0.
\end{equation}
That is, for constant, $\rho$, $T$ and $\theta$, the complementary free energy  tends to a maximum in equilibrium.

To make our approach more clear we  ignore the thermal variables and focus on the purely mechanical case, that is we assume that the viscoelastic material undergoes an isothermal process in which \eqref{CD-gibbs} reduces to
\begin{equation}\label{CD-isoth}
    \rho \dot{\phi_{c}} + TA\frac{\dot{\rho}}{\rho} -\dot{T}A \geq 0.
\end{equation}
In what follows we will focus on  providing a thermodynamic framework using the complementary free energy (or the Gibbs free energy). In particular,  we will show that \eqref{CD-isoth}  places restrictions on the constitutive relation for the strain $A$.

In order to explain the response of rate-type viscoelastic solids, we now assume that the  set of variables $\phi_{c}$, $A$ are functions of the set of independent variables $\rho$, $T$, $\dot{T}$:
\begin{equation}
  \phi_{c}= \phi_{c}(\rho, T, \dot{T}),~~~~   A = A(\rho, T, \dot{T}).
 \end{equation}
We also assume that $\phi_{c}$ and $A$ are differentiable in $\rho$, $T$ and $\dot{T}$ as many times as we need. We further suppose that the complementary free energy and the strain are identically zero when both the stress and the stress rate are zero  and the density is the reference density $\rho_{0}$ (a stress and stress-rate free state). That is, we assume that $\phi_{c}(\rho_{0},0,0)\equiv 0$ and $A(\rho_{0},0, 0)\equiv 0$.  Now, if we calculate the material time derivative of $\phi_{c}$ and substitute the resulting expression into \eqref{CD-isoth}, we get
\begin{equation}\label{GCD}
    \frac{1}{\rho}\left(\rho^{2}\frac{\partial \phi_{c}}{\partial \rho}+TA\right)\dot{\rho}
        +\left(\rho \frac{\partial \phi_{c}}{\partial T} -A \right) \dot{T}
        + \rho  \frac{\partial \phi_{c}}{\partial \dot{T}}\ddot{T} \geq 0.
\end{equation}
Assuming that $\dot{\rho}$ and $\ddot{T}$ may take arbitrary values and requiring that  the above condition needs to hold for all $\dot{\rho}$ and $\ddot{T}$, we see that the coefficients of $\dot{\rho}$ and $\ddot{T}$  must vanish. This leads to the constraints
\begin{equation}\label{constraint-G}
   \rho^{2}\frac{\partial \phi_{c}}{\partial \rho}+TA=0, ~~~~ \frac{\partial \phi_{c}}{\partial \dot{T}}=0
\end{equation}
and \eqref{GCD} becomes
\begin{equation}\label{redCD}
        \left(\rho \frac{\partial \phi_{c}}{\partial T} -A \right) \dot{T}\geq 0.
\end{equation}
In view of the second constraint in \eqref{constraint-G}, the response function for the complementary free energy  is independent of the stress-rate: $\phi_{c}= \phi_{c}(\rho, T)$. If we denote the coefficient of $\dot{T}$ in \eqref{redCD}  by $A_{d}$, it  becomes
\begin{equation}\label{diss-ineq}
    A_{d}(\rho, T, \dot{T}) \dot{T} \geq 0.
\end{equation}
Then the constitutive equation for $A$ can be expressed as
\begin{equation}\label{constitutive}
    A(\rho, T, \dot{T})=A_{e}(\rho, T)- A_{d}(\rho, T, \dot{T}),~~~~A_{e}(\rho, T)=\rho \frac{\partial \phi_{c}(\rho, T)}{\partial T},
\end{equation}
where $A_{e}$ and $A_{d}$ stand for the elastic and dissipative parts of the strain and they are measures of the elastic and viscoelastic responses of the material. We  require that  the entropy production is the minimum, namely zero, at the states such that $\dot{T} = 0$. So we set $A_{d}(\rho, T, 0) = 0$ and $A(\rho, T, 0)=A_{e}(\rho, T)$.

Recalling that the relation between the Helmholtz free energy $\psi$ and the Gibbs free energy  $G$  (per unit mass) is given by  $\psi = G - T \frac{\partial G}{\partial T}$, we obtain from \eqref{comp} the relation between $\phi_{c}$ and $G$ as
\begin{equation}\label{gibbs}
  \phi_{c}= \frac{TA}{\rho}- G + T \frac{\partial G}{\partial T}.
 \end{equation}
 It is worth to note that we could also write the Clasius-Duhem inequality \eqref{CD-isoth} in terms of $G(\rho, T, \dot{T})$. Using \eqref{constitutive} and \eqref{gibbs}, we can state $A_{d}$ in the form
\begin{equation}\label{compl}
    A_{d}= T \left(\frac{\partial A}{\partial T}-C\right),~~~~ C=-\rho  \frac{\partial^{2} G}{\partial T^{2}},
\end{equation}
and the inequality  \eqref{diss-ineq} becomes $T\left(\frac{\partial A}{\partial T}-C\right)\dot{T}\geq 0$. We shall later show how the function $C$ is related to the instantaneous elastic compliance. Using  \eqref{constitutive} in the first constraint of \eqref{constraint-G} we can restate it as
\begin{equation}\label{pressure}
  (T+P)\rho\frac{\partial \phi_{c}}{\partial T}- TA_{d}=0, ~~~~
  P(\rho, T)=\rho \frac{\partial \phi_{c}/\partial \rho}{\partial \phi_{c}/\partial T},
\end{equation}
where the function $P$ is the thermodynamic pressure of the viscoelastic material.

We now restrict ourselves to a particular case in which the strain $A$ has a dissipative part that is linear in $\dot{T}$. From an analytical point of view, the simplest form is that $A_{d}(\rho, T, \dot{T}) = \gamma(\rho, T) \,\dot{T}$ with $\gamma$ being a scalar function. As a result, we have
\begin{equation}\label{new-model-dot}
    A(\rho, T, \dot{T}) = h(\rho, T) - \gamma(\rho, T) \,\dot{T},
\end{equation}
where we use the notation $h(\rho, T)=\rho \partial \phi_{c}(\rho, T)/\partial T$.  Note that, with this choice of $A_{d}$, the necessary and sufficient condition that the thermodynamic inequality \eqref{diss-ineq} is valid is $\gamma(\rho, T) \geq 0$. Since $A(\rho_{0},0, 0)= 0$ and $A_{d}(\rho_{0}, 0, 0) = 0$ in the absence of the stress and the stress-rate, we require that $h(\rho_{0},0) = 0$.

\subsection{Strain-limiting viscoelastic model}

As in the previous sections we now require that  $\epsilon=\partial u/\partial x$ is small for all $x$ and $t$ considered. This allows us to replace $A$ by $\epsilon$ in the equations derived above. Noting that  the balance of mass  given by $\rho (1+A)=\rho_{0}$ in one-dimensional case reduces to $\rho(1+\epsilon)=\rho_{0}$ due to the smallness assumption,  we reach the conclusion that $\rho$ and $\rho_{0}$ are of the same order if $\epsilon$ is small enough. So, in the context of strain-limiting theory, one may use $\rho$ and $\rho_{0}$ interchangeably. Furthermore, with the small strain assumption, the difference between the quantities measured in the reference and current configurations disappears and the material time derivatives become partial derivatives with respect to time. Consequently, in the strain-limiting case,   \eqref{new-model-dot} reduces to \eqref{new-model} with $\gamma \geq 0$.

\subsection{Strain-limiting elastic model}
Proceeding similarly, we could also derive the constitutive relation for strain-limiting elastic materials. In the case  of an elastic material the response functions for $\phi_{c}$ and $A$ depend only on the density $\rho$ and the stress $T$  and are independent of their rates. So, assuming  $\phi_{c} = \phi_{c}(\rho,T)$ and $A = A(\rho,T)$ (and consequently $G = G(\rho,T)$) and using \eqref{CD-isoth}, we get
\begin{displaymath}
    \frac{1}{\rho}\left(\rho^{2}\frac{\partial \phi_{c}}{\partial \rho}+TA\right)\dot{\rho}
        +\left(\rho \frac{\partial \phi_{c}}{\partial T} -A \right) \dot{T}
         \geq 0
\end{displaymath}
instead of \eqref{GCD}. Since  $\dot{T}$ is not  the independent variable of $\phi_{c}$ and $A$, we have
\begin{equation}\label{elas}
    A(\rho, T)=\rho \frac{\partial \phi_{c}(\rho, T)}{\partial T},
\end{equation}
(that is, $A\equiv A_{e}$ and $A_{d}\equiv 0$) which delivers the constitutive equation for the strain-limiting model of elastic materials.  Since $\dot{\rho}$ is not  the independent variable of $\phi_{c}$ and $A$, we have the  first constraint in \eqref{constraint-G}. A simple calculation shows that it is satisfied identically if the thermodynamic pressure is given by $P=-T$. In the present case the quantity $C$ defined in \eqref{compl} becomes the instantaneous elastic compliance. Furthermore, using \eqref{gibbs} in \eqref{elas}, we observe that, for elastic materials, the complementary free energy is the negative of the Gibbs energy: $\phi_{c}=-G$. Note that the above constitutive relation \eqref{elas} for $A$ can be also obtained as the limiting case of the constitutive relation \eqref{constitutive} for viscoelastic solids. That is, the elastic behavior of the model can be simply recovered through assuming that  $\dot{T} \equiv 0$ in \eqref{constitutive}.

For strain-limiting elastic solids, approximating $A$ by $\epsilon$ in \eqref{elas} we get $\epsilon = h(T)$, where $h$ is as before. This is in full agreement with  \cite{Rajagopal2011a}. It is worth noting that, in the elastic case,  the  constitutive relation for $\epsilon$ is based on the existence of a scalar function $\phi_{c}$ depending on  $T$. A similar result  was already proposed in \cite{Bustamante2011}.

\section{Comparison of two strain-limiting models of viscoelasticity}\label{sec:sec5}

In this section we will compare the strain-rate type model \eqref{Raj-model} and the stress-rate type model \eqref{new-model} of strain-limiting viscoelastic materials. In the rest of this section we will compare the models \eqref{Raj-model} and  \eqref{new-model} from various aspects: energy decay, the nonlinear differential equations of motion and Fourier analysis of the corresponding linear models.

\subsection{Energy decay}

 In one space dimension, the equation of motion for a homogeneous, viscoelastic, infinite medium  is given by
\begin{equation}\label{eqn-motion}
    \rho\, u_{tt} = T_{x}.
\end{equation}
We assume that the density is constant and the stress tends to zero at infinity: $T \rightarrow 0$ as $x \rightarrow \pm \infty$. Below we show that  the energy decays over time for both the strain-rate type model \eqref{Raj-model} and the stress-rate type model \eqref{new-model}.

We start with the stress-rate type model \eqref{new-model}. Multiplying \eqref{eqn-motion} by $u_{t}$ and integrating over the space variable, we get
\begin{equation}
    \frac{1}{2}\frac{d}{dt} \int_{-\infty}^{\infty} \rho (u_{t})^{2}\,dx = \int_{-\infty}^{\infty} T_{x} u_{t} \,dx.
\end{equation}
One integration by parts on the  right-hand side of this equation yields
\begin{equation}\label{stress-energy}
    \frac{1}{2} \frac{d}{dt} \int_{-\infty}^{\infty} \rho (u_{t})^{2}\,dx     = - \int_{-\infty}^{\infty}T \epsilon_{t}  \,dx.
\end{equation}
Rewriting the integrand on the right-hand side as $T \epsilon_{t}=(T \epsilon)_{t}-  T_{t} \epsilon  $ and then using the constitutive relation \eqref{new-model} in this equation we get $T \epsilon_{t}=(T \epsilon)_{t}-h(T)T_{t}+\gamma (T_{t})^{2}$. Recalling that $h(T)=\rho  d\phi_{c}(T)/d T$ and noticing that $h(T)T_{t}=  (\rho\phi_{c}(T))_{t}$ we see this equation is equivalent to $T \epsilon_{t}=(T \epsilon-\rho  \phi_{c})_{t}+\gamma (T_{t})^{2}$. Finally, to express it in terms of the internal energy $\omega(T)$ (or equivalently in terms of the Helmholtz free energy $\psi(T)$) we use the linearized form $\rho\phi_{c}=T\epsilon-\rho \psi$ of \eqref{comp} and we reach $T \epsilon_{t}=\rho  \omega_{t}+\gamma (T_{t})^{2}$. Thus \eqref{stress-energy} takes the form
\begin{displaymath}
    \frac{d}{dt} \int_{-\infty}^{\infty} \left(\frac{\rho}{2} (u_{t})^{2} + \rho\, \omega(T) \right)  dx = - \gamma \int_{-\infty}^{\infty} (T_{t})^{2} \,dx,
\end{displaymath}
which, since the right-hand side is always negative for $\gamma>0$, shows that the total energy given by the integral on the left-hand side is decaying over time.

For the strain-rate type model \eqref{Raj-model} we proceed similarly to get \eqref{stress-energy}. Then, rewriting  the equation of motion \eqref{eqn-motion} in the form $\rho\, \epsilon_{tt} = T_{xx}$  in terms of the strain and using the constitutive relation \eqref{Raj-model}, one can indeed verify that $\epsilon_{t}=(g(T))_{t}-(\nu/\rho)T_{xx}$ and consequently  $T\epsilon_{t}=T(g(T))_{t}-(\nu/\rho)TT_{xx}$. Assuming the existence of a potential function $\bar{G}(T)$ with $g(T)=-\rho  d\bar{G}(T)/d T$ we get  $T\epsilon_{t}=\rho (\bar{G}-Td\bar{G}(T)/dT)_{t}-(\nu/\rho)TT_{xx}$. This can be written as $T\epsilon_{t}=\rho \omega_{t}-(\nu/\rho)TT_{xx}$   in terms of the internal energy $\omega(T)$ (or equivalently in terms of the Helmholtz free energy $\psi(T)$) if we require that the potential function $\bar{G}(T)$ is the Gibbs free energy $G(T)$. Substitution of the above relation obtained for  $T\epsilon_{t}$ into \eqref{stress-energy}  yields
\begin{displaymath}
    \frac{d}{dt} \int_{-\infty}^{\infty} \left(\frac{\rho}{2} (u_{t})^{2} + \rho\, \omega(T) \right)  dx = - \frac{\nu}{\rho} \int_{-\infty}^{\infty} (T_{x})^{2} \,dx,
\end{displaymath}
where one integration by parts has been performed on the right-hand side. For $\nu >0$, this result shows again  that the total energy  is decreasing over time.

\subsection{Nonlinear differential equations of motion}

We now define the dimensionless quantities
\begin{equation}\label{dim}
    \bar{x} = \frac{x}{L},\,\,\, \bar{t} = \frac{t}{L}\sqrt{\frac{\mu}{\rho}},\,\,\, \bar{T}= \frac{T}{\mu},  \,\,\,
        \bar{\nu} = \frac{\nu}{L} \sqrt{\frac{\mu}{\rho}},\,\,\,\bar{\gamma} = \frac{\gamma \mu}{L} \sqrt{\frac{\mu}{\rho}},
 \end{equation}
where $L$ is a characteristic length and $\mu$ is a constant with the dimension of the stress.  Substituting \eqref{new-model} into the equation of motion for the strain, $\rho \epsilon_{tt} = T_{xx}$,  and using \eqref{dim}  in the resulting equation  we get
\begin{equation}\label{new-ode}
    T_{xx} + \gamma\,T_{ttt} = h(T)_{tt},
\end{equation}
where we have dropped the overbars for simplicity. On the other hand, following a similar approach, the model \eqref{Raj-model} gives
\begin{equation}\label{old-ode}
    T_{xx} + \nu T_{xxt} = g(T)_{tt}
\end{equation}
(see e.g. (1.1) in \cite{Erbay2015} for details).  The above two differential equations, \eqref{new-ode} and \eqref{old-ode}, are the main equations that characterize the behaviour of viscoelastic solids for two different strain-limiting models. Both equations differ from the corresponding equations within the context of classical viscoelasticity theories in the sense that the they are written in terms of the stress and also the non-linearity is on the inertia term.

We now consider the travelling wave solutions of \eqref{new-ode} and \eqref{old-ode}. If we substitute $T = T(\xi)$ with $\xi = x - c t,$ where $c$ is the wave propagation speed, \eqref{new-ode} and \eqref{old-ode} both give
\begin{equation}\label{ode-tw}
    T^{\prime\prime} -  \kappa T^{\prime\prime\prime} = c^{2} [f(T)]^{\prime\prime},
\end{equation}
(with $f=h$, $\kappa = \gamma c^{3}$ for \eqref{new-ode} and  with  $f=g$,  $\kappa = \nu c$ for \eqref{old-ode}) where the symbol $^{\prime}$ denotes ordinary differentiation. We note that \eqref{ode-tw}  is exactly the same as equation $(3.2)$ of \cite{Erbay2015}. Erbay and \c{S}eng\"ul \cite{Erbay2015} investigated existence of travelling wave solutions of  \eqref{ode-tw} corresponding to the heteroclinic connections between two constant equilibrium states at infinity. We arrive at the conclusion that all observations made in  \cite{Erbay2015} for the model \eqref{Raj-model} are also valid for the new model \eqref{new-model}.

\subsection{Fourier analysis of the corresponding linear models}

If we assume that the functions $h$ and $g$ in  \eqref{new-model} and \eqref{Raj-model}, respectively, are linear functions of the stress, that is, $h(T) = h^{\prime}(0) T$ and $g(T) = g^{\prime}(0) T$, equations \eqref{new-ode} and \eqref{old-ode} become
\begin{equation}\label{new-ode-linear}
    T_{xx} + \gamma\,T_{ttt} = h^{\prime}(0) T_{tt}
\end{equation}
and
\begin{equation}\label{old-ode-linear}
 T_{xx} + \nu T_{xxt} = g^{\prime}(0) T_{tt},
\end{equation}
 respectively. For simplicity, throughout the rest of this work we take $h^{\prime}(0)=g^{\prime}(0)  = 1$.

In order to investigate the Hadamard stability/instability, we now do a Fourier mode analysis for both  \eqref{new-ode-linear} and \eqref{old-ode-linear}. To do this, we consider a spatial mode with a wavenumber $k$, temporal growth rate $r$, and amplitude $T_{a}$ so that $T(x, t)  = \mathrm{Re}\{ T_{a} e^{r t} e^{i k x} \}$.

Considering  \eqref{old-ode-linear} first, we obtain $ r^{2} + \nu k^{2} r + k^{2} = 0$ which is a quadratic in $r$, and hence
 \begin{displaymath}
    r = \frac{1}{2}(- \nu k^{2} \mp \sqrt{\Delta}),
\end{displaymath}
where $\Delta = k^{2}( \nu^{2} k^{2} - 4 )$. Clearly, there exists a critical value for $k$, which can be found as $k_{c}^{2} = 4/ \nu^{2}$. Therefore, if $k^{2} < k_{c}^{2}$ we have $\Delta < 0$ and there exists two complex conjugate roots satisfying $\mathrm{Re}(r)  = - \nu k^{2}/{2} < 0$ leading to stability. On the other hand, if $k^{2} > k_{c}^{2}$ we have $\Delta > 0$ and we obtain two real roots $r_{1} $ and $r_{2} $. Since they are both negative we have stability again.

Considering \eqref{new-ode-linear} now, we obtain the cubic equation $\gamma r^{3} - r^{2} - k^{2} = 0$. Since the discriminant for this cubic equation is $- k^{2} (4 + 27 \gamma^{2} k^{2})$, which is always negative, we can immediately say that the cubic equation has one real root and two complex conjugate roots. Without calculating the roots, we can apply Descartes' rule of signs for polynomials to say that the only real root is positive and therefore there is an instability in this case. As a result, while the solutions of \eqref{old-ode-linear} are stable for any wave number $k$, the solutions of \eqref{new-ode-linear} are unstable for all $k$.


\end{document}